\documentclass{amsart}
\usepackage{amssymb,latexsym,amsmath,amsfonts}
\newtheorem{definition}{Definition}[section]
\newtheorem{theorem}[definition]{Theorem}
\newtheorem{corol}[definition]{Corollary}
\newtheorem{lemma}[definition]{Lemma}

\newcommand{\F}{\mathbb{F}}
\newcommand{\fq}{\mathbb{F}_q}
\newcommand{\rmv}[1]{}

\begin{document}

\title{On explicit factors of Cyclotomic polynomials over finite fields}
\author{Liping Wang}
\address{Center for Advanced Study,
Tsinghua University, HaiDian District,
Beijing(100084), China.}
\email{wanglp@mail.tsinghua.edu.cn}

\author{Qiang Wang}
\address{School of Mathematics and Statistics, Carleton
University, 1125 Colonel By Drive, Ottawa, Ontario, K1S 5B6,
Canada.} 
\email{wang@math.carleton.ca} 
\thanks{Research is partially supported by NSERC of Canada. }

\begin{abstract}
We study the explicit factorization of $2^n r$-th cyclotomic polynomials over finite field $\fq$ where $q, r$ are odd
 with $(r, q) =1$. We show that all irreducible factors of  $2^n r$-th cyclotomic polynomials can be obtained easily from irreducible factors of cyclotomic polynomials of small orders. In particular, we obtain the explicit factorization of $2^n 5$-th cyclotomic polynomials over finite fields and construct several classes of irreducible polynomials of degree $2^{n-2}$  with fewer than $5$ terms.  The reciprocals of these irreducible polynomials are irreducible polynomials of the form $x^{2^{n-2}} + g(x)$ such that the degree of $g(x)$ is small ($\leq 4$), which could have potential applications as mentioned by Gao, Howell, and Panario in \cite{GaoHowellPanario}.
\end{abstract}

\maketitle

\section{Introduction}
Let $p$ be prime, $q=p^m$, and $\F_q$ be a finite field of order
$q$.  Let $Q_{n}(x)$ denote the $n$-th cyclotomic polynomial 
\[
Q_{n}(x)= \prod_{0<j \leq n,  (j,n)=1}(x- \zeta^j) 
\]
where $\zeta$ is a primitive $n$-th root of unity. Clearly $x^n - 1 = \prod_{d|n} Q_d(x)$
and the M\"{o}bius inversion formula  gives $Q_n(x) = \prod_{d|n} (x^d - 1)^{\mu(n/d)}$ where $\mu$ is the M\"{o}bius function.  If $(q, n) =1$, then it is well known that $Q_n(x)$ can be factorized into $\phi(n)/d$ distinct monic irreducible polynomials of the same degree $d$ over $\fq$, where $d$ is the least positive integer such that $q^d \equiv 1 \pmod{n}$ (see \cite[Theorem 2.47]{LN}).  Basically we know the number and the degree of irreducible factors of cyclotomic polynomials. However, factoring cyclotomic polynomials $Q_n(x)$ over finite field $\fq$  explicitly  still remains as a fundamental question.  Moreover, it is also known that  explicit factorization of cyclotomic polynomials is related to  the factorization of   other interesting classes of  polynomials.  For example, Fitzgerald and Yucas \cite{FY1} have discovered a nice link between  the factors of Dickson polynomials over finite fields  and factors of cyclotomic polynomials and self-reciprocal polynomials recently.  This means that factoring cyclotomic polynomials explicitly provides an alternative way to factor Dickson polynomials explicitly. 

 Explicit factorization of $2^n$-th cyclotomic polynomials  $Q_{2^n}(x)$ over $\fq$ are given in \cite{LN}
 when $q\equiv 1 \pmod{4}$ and in \cite{Meyn} when $q\equiv 3 \pmod{4}$.   Recently, Fitzgerald and Yucas   
\cite{FY2} have studied explicit factors of  $2^n r$-th cyclotomic polynomials $Q_{2^nr}(x)$  where $r$ is 
prime and $q\equiv \pm 1 \pmod{r}$  
over finite field $\fq$ in order to obtain explicit factorization of Dickson polynomials.
This  gives a complete answer to the explicit  factorization of cyclotomic polynomials $Q_{2^n3}(x)$ and thus  Dickson polynomials  $D_{2^n 3}(x)$ of the first kind over $\fq$. 
However, the general situation for arbitrary  $r$ remains open.  
Without loss of generality  we  assume that $(2r, q) =1$. In this paper, we reduce the problem of factorizing all $2^n r$-th cyclotomic 
polynomials over $\fq$ into factorizing a finite number of lower degree cyclotomic polynomials over $\fq$.  In particular, 
we  give the explicit factorization of cyclotomic polynomials $Q_{2^n r}(x)$  
over $\fq$ where $r=5$.   
The method we are using is a combination of case analysis, factorizing low degree polynomials, and
the recursive construction based on basic properties of cyclotomic polynomials.

The irreducible factors of these cyclotomic polynomials are sparse polynomials (polynomials with a few nonzero terms).  Sparse irreducible polynomials 
 are important in efficient hardware implementation of feedback shift registers and finite field arithmetic (\cite{Berlekamp}, \cite{GolombGong},
\cite{WangBlake}). The second focus of our paper is to explicitly construct sparse irreducible polynomials of high degrees. We remark that explicit construction
of irreducible polynomials in general has attracted a lot of attentions and a lot of progress has made in the past two decades. Most of these
constructions are iterated constructions which extend the classical transformation $f(x) \rightarrow f(x^n)$. A
a nice survey  on this topic as of year 2005 can be found in \cite{Cohen}. Here we are interested in sparse irreducible 
polynomials and the classical transformation is used. Therefore the main tool in the paper 
is the following classical result which help us to   
construct high degree irreducible polynomials based on low degree irreducible polynomials.

\begin{lemma}[Theorem 3.35 in \cite{LN}]
\label{irred}
Let $f_1(x)$, $f_2(x)$, \ldots, $f_N(x)$ be all distinct monic irreducible polynomials in $\fq[x]$ of degree $m$ and order $e$, and let $t \geq 2$ be an integer whose prime factors divide $e$ but not $(q^m-1)/e$. Assume also that $q^m \equiv 1 \pmod{4}$ if $t \equiv 0 \pmod{4}$. Then $f_1(x^t), f_2(x^t), \ldots, f_N(x^t)$ are all distinct monic irreducible polynomials in $\fq[x]$ of degree $mt$ and order $et$.
\end{lemma}

In Section~\ref{methodology}, we describe the methodology used in this paper to factor $Q_{2^n r}(x)$ over $\fq$. We prove that all irreducible factors of $Q_{2^n r}(x)$ can be obtained easily from irreducible factors of 
$Q_{2^L r}(x)$ where $L$  is a small constant depending on $q$ and $r$ (Theorem~\ref{general}). This also provide us a way 
to construct sparse irreducible polynomials of high degree $2^n r$ over $\fq$. We note that the result in this section is 
true for any odd $q, r$ such that $(q, r)=1$. Then the rest of paper deals with $r=5$. 
In Section~\ref{old}, we obtained the factorization results of $Q_{2^n 5}(x)$ when $q \equiv \pm 1 \pmod{5}$. These results are not new and can also be found in \cite{FY2}. However, for the sake of completeness, we also include the proofs.  As a consequence,  we obtain several
classes of  irreducible binomials/trinomials of degree $2^{n-2}$ over $\fq$ where $q \equiv \pm 1 \pmod{5}$.  
In Section~\ref{new1}, we consider the situation when $q\equiv 13, 17 \pmod{20}$.  We obtain the explicit factorization of $Q_{2^n 5}(x)$  
in Theorems~\ref{factor-mod-13-17}. Moreover, we can construct several classes of irreducible polynomial of 
five terms with degree $2^m$ 
(Corollary~\ref{irred-mod-13-17}). Then the case  of $q\equiv 3 \pmod{20}$  is considered in Section~\ref{new2}. 
The factorization results are given in Theorem~\ref{factor-mod-3} for 
$q\equiv 3 \pmod{20}$ such that $q\neq 3$ and 
in Theorem~\ref{factor-3} for $q=3$.  The irreducible polynomials constructed are trinomials and pentanomials (see Corollaries~\ref{irred-mod-3},~\ref{irred-3}).  In Section~\ref{new3}, the case of $q \equiv 7 \pmod{20}$ is considered  and the result can be found in Theorem~\ref{factor-mod-7}. 

 We note that the reciprocals of sparse irreducible polynomials constructed in this paper can be written as the form of $x^n + g(x)$ where the degree of $g(x)$ is at most $4$.  It is well known that irreducible polynomials of form $x^n + g(x)$ with $g(x)$ having a small degree are desirable in implementing pseudorandom number generators and in constructing elements of provable high orders in finite fields (see the survey paper \cite{GaoHowellPanario}).    Therefore our irreducible polynomials might be useful in some of  applications  mentioned in \cite{GaoHowellPanario}.

\section{Methodology and notations}\label{methodology}

In this section, we  describe the method that we are using in this paper.  First of all we recall the following basic 
results on cyclotomic polynomials.
 
\begin{lemma}\cite[Exercise~2.57]{LN}
\label{cyclo}

(a) $Q_{2n}(x) = Q_n(-x)$ for $n \geq 3$ and $n$ odd.

(b) $Q_{mt}(x) = Q_{m}(x^t)$ for all positive integers $m$ that are divisible by the prime $t$.

(c) $Q_{mt^k}(x) = Q_{mt}(x^{t^{k-1}})$ if $t$ is a prime and $m, k$ are arbitrary positive integers. 
\end{lemma}

Let us start with the 
factorizations of $Q_r(x)$ and $Q_{2r}(x) = Q_r(-x)$.  
 Because of Lemma~\ref{cyclo}, we have $Q_{2^n r} (x) = Q_{2^{n-1} r}(x^2)$ for $n \geq 2$.  
Hence the key to continue the process of factorization is to factor $Q_{2^{n-1} r}(x^2)$ into a product of irreducible polynomials once we obtain the factorization of $Q_{2^{n-1} r}(x)$.

Now, we show that we can reach to the end after only a finite number of iterations. 
 Let $v_2(k)$ denotes the highest power of $2$ dividing $k$ and $L_i = v_2(q^i -1)$ for $ i \geq 1$. 
In particular, let $L:= L_{\phi(r)}= v_2(q^{\phi(r)}-1)$ where $\phi$ is the Euler's phi function.
 Then we have the following result.

\begin{theorem}\label{general}
Let   $q=p^m$ be a power of an odd prime $p$,
let $r \geq 3$ be  any odd number such that $(r, q) =1$,  and let
$L := L_{\phi(r)}= v_2(q^{\phi(r)}-1)$, the highest power of $2$ dividing $q^{\phi(r)} -1$ with $\phi(r)$ the Euler's phi function.
For any $n\geq L$ and any  irreducible factor $f(x)$ of  $Q_{2^L r}(x)$ over $\fq$,  $f(x^{2^{n-L}})$ is also 
irreducible over $\fq$. Moreover, all irreducible factors of  $Q_{2^n r}(x)$ are obtained  in this way.
\end{theorem}

\begin{proof}
Because $q, r$ are odd, we have $\phi(r)$ is even and then  $n \geq L \geq 2$. By \cite[Theorem 2.47]{LN},  
$2^{L} r$-th cyclotomic polynomial $Q_{2^L r}(x)$ has $\phi(2^{L}  r)/m$ distinct monic irreducible factors of the same degree $m$,
 where $m$ is the least positive integer  such that $q^m \equiv 1 \pmod{2^{L}  r}$. 
Because $q^{\phi(r)} \equiv 1 \pmod{r}$ and $L = v_2(q^{\phi(r)}-1)$,  we have
$q^{\phi(r)} \equiv 1 \pmod{2^L r}$. This implies that $m \leq \phi(r)$.  By the definition of $L$, we must have
 $2^{L+1} \nmid (q^m -1)$. Since 
each factor has order $e= 2^{L}  r$ and $2 \nmid (q^m -1)/e$,  by Lemma~\ref{irred},  each 
irreducible polynomial $f(x)$ of $Q_{2^L r}(x)$  generates an irreducible factor $f(x^2)$ of $Q_{2^{L+1} r}(x)$. 
More generally, $f(x^{2^{n-L}})$ is also irreducible factor of $Q_{2^n r}(x)$  since $L \geq 2$ implies 
that $4 \mid q^{m} -1$.  Moreover, $f(x^{2^{n-L}})$ has degree $m2^{(n-L)}$ and order $2^{n} r$. Hence there are $\phi(2^n r)/(m2^{n-L}) = 2^{n-1} \phi(r) / (m2^{n-L}) = \phi(2^L r)/m$ distinct 
irreducible factors for $Q_{2^n r}(x)$. Therefore all irreducible factors of  $Q_{2^n r}(x)$  are constructed from irreducible factors of $Q_{2^L r}(x)$ over $\fq$.
\end{proof}

Theorem~\ref{general} tells us that a recursive way of factoring  $2^n r$-th cyclotomic polynomials essentially requires only finitely many factorizations of low degree polynomials (at most $L$ iterations starting from $Q_{r}(x)$). This also provides us a method to construct irreducible polynomials from low degree irreducible polynomials. The fact that we use the classical transformation on low degree polynomials can guarantee  the resulting high degree irreducible polynomials are sparse polynomials. 

For $n < L$,  since each irreducible factor $f(x)$ of $Q_{2^{n-1} r}(x)$ has the same degree $m$ and $f(x^2)$ may not be irreducible polynomial of degree $2m$, we need to factor $f(x^2)$ further. And in most cases, we need to factor $f(x^2)$ into two irreducible polynomials of degree $m$.  
We  see more in detail for $r=5$ in the forthcoming sections. 
This involves the process of factoring certain types of polynomials of degree $8$ into two quartic polynomials.

Finally we fix some other notations for the rest of the paper.

Let $\Omega(k)$ denote the set of primitive $k$-th root of unity.  In particular, $\Omega(2^0) = \{ 1\}$, $\Omega(2^1) = \{-1\}$.
Let  $\rho_n$ denote an arbitrary element in $\Omega(2^n)$.

The expression $\prod_{a \in  A} \cdots \prod_{b \in B}  f_i(x, a, \ldots, b)$ denotes the product of distinct
 irreducible polynomials $f_i(x, a, \ldots, b)$ satisfying conditions $a\in A$, $\ldots$,  $b \in B$.  

Let $\left( a \atop p \right)$ denote the Legendre symbol and the following basic results on Legendre symbols are also used in the paper. 

(i) $\left( 2 \atop p \right) = 1$ if and only if $p \equiv 1, 7 \pmod 8$;

(ii)  $\left( -2 \atop p \right) = 1$ if and only if $p \equiv 1, 3 \pmod 8$;

(iii) $\left( 5 \atop p \right) = 1$ if and only if $p \equiv \pm 1, \pm 9 \pmod{20}$. 

\section{Case: $q\equiv \pm 1 \pmod{5}$}
\label{old}

Recall that $L_i = v_2(q^i -1)$, the highest power of $2$ dividing $q^i -1$ for $ i \geq 1$. If $q\equiv \pm 1 \pmod 5$ and $q\equiv 1 \pmod 4$ (i.e., $q \equiv 1 \pmod{20}$ or $q\equiv 9 \pmod{20}$), then $L_4 = L_2+1$ and $L_2 = L_1 + 1$. Moreover, $\rho_1 = -1$ must be a square and thus $\rho_2^2 = \rho_1$. 

Similarly, if $q\equiv \pm 1 \pmod 5$ and $q\equiv 3 \pmod 4$ (i.e., $q \equiv 11 \pmod{20}$ or $q\equiv 19 \pmod{20}$), then $L_4 = L_2+2$ and $L_2 = L_1 + 1$. Moreover, $\rho_1 = -1$ can not be a square.

We have the following results for these four different cases. 

\begin{theorem}
Let $q\equiv  1 \pmod{20}$. 
Then we have the following factorization of $2^n5$-th cyclotomic polynomial $Q_{2^n 5} (x)$  over $\fq$.  

(i) $Q_{5}(x) = \prod_{w \in \Omega(5)} (x-w)$ and $Q_{10}(x) =    \prod_{w \in \Omega(5)} (x+ w)$.

(ii) If $2 \leq n \leq L_1$, then 
\[
Q_{2^n 5} (x) = \prod_{w \in \Omega(5)} \prod_{\rho_n \in \Omega(2^n)}  \left( x  - w  \rho_n \right).
\]

(iii)if $n \geq L_2=L_1+1$, then 
\[
Q_{2^n 5} (x) =
\prod_{w \in \Omega(5)} \prod_{\rho_{L_1} \in \Omega(2^{L_1})}   \left( x^{2^{n-L_1}} - w  \rho_{L_1} \right).
\]
\rmv{

(ii) if $n = L_2 = L_1 +1$, then 
\[
Q_{2^n 5} (x) = \prod_{w \in \Omega(5)} \prod_{\rho_{L_1} \in \Omega(2^{L_1})}  \left( x^2 - w \rho_{L_1} \right).
\]

(iii) if $n = L_4 = L_2 +1$, then 
\[
Q_{2^n 5} (x) = \prod_{w \in \Omega(5)} \prod_{\rho_{L_1} \in \Omega(2^{L_1})}  \left( x^4 - w \rho_{L_1} \right).
\]

(iv)  if $n > L_4 = L$, then 

\[
Q_{2^n 5} (x) = Q_{2^L 5}(x^{2^{n-L}}) = \prod_{w \in \Omega(5)} \prod_{\rho_{L_1} \in \Omega(2^{L_1})}  \left( x^{2^{n-L+2}} - w \rho_{L_1} \right).
\]

}
\end{theorem}

\begin{proof}
In this case, $5 \mid q-1$. Hence $\Omega(5) \subseteq \fq$. Moreover,  $\Omega(2^n) \subseteq \fq$ for all $n \leq L_1$. Moreover, if $n \leq L_1$ and $\rho_{n-1} \in \Omega(2^{n-1})$ then we can find $\rho_{n} \in  \Omega(2^n) \subseteq \fq$ such that $\rho_n^2 = \rho_{n-1}$. 

(i) Because $\rho_0 =1$ and $\rho_1 = -1$,  it is obvious
to see that  $Q_5 (x) = \prod_{w \in \Omega(5)}   \left( x - w \right) =\prod_{w \in \Omega(5)}  \left( x - w \rho_0 \right)$ and
$Q_{10}(x) = Q_5(-x) = \prod_{w \in \Omega(5)}  \left( x + w \rho_0 \right)= \prod_{w \in \Omega(5)}   \left( x - w \rho_1 \right)$.

(ii) For $ 2\leq n \leq L_1$, we have
$Q_{2^n 5}(x) = Q_{2^{n-1} 5}(x^2) = \prod_{w \in \Omega(5)} \prod_{\rho_{n-1} \in \Omega(2^{n-1})}  
\left( x^2 - w \rho_{n-1} \right)$ by induction hypothesis and Lemma~\ref{cyclo}. 
Because each $w\in \Omega(5)$ can be written as $w=u^2$ where $u\in \Omega(5)$, we obtain that
$x^2 - w \rho_{n-1} = x^2 - u^2 \rho_{n}^2 = (x-u \rho_{n})(x+ u \rho_{n})$ and both $\rho_{n}, -\rho_{n} \in \Omega(2^n)$. Hence,  for $ 2\leq n \leq L_1$, we have
\[
Q_{2^n 5}(x)
= \prod_{u \in \Omega(5)} \prod_{\rho_{n} \in \Omega(2^{n})}  
 \left( x - u \rho_{n} \right). 
\]

(iii) Again, we have
$Q_{2^{L_2} 5}(x) = Q_{2^{L_1} 5}(x^2) = \prod_{w \in \Omega(5)} \prod_{\rho_{L_1} \in \Omega(2^{L_1})}   \left( x^2  -  w \rho_{L_1} \right)$. 
Because $w$ is a square element and $\rho_{L_1}$ is a non-square element in $\fq$, $x^2 - w  \rho_{L_1}$ is an irreducible polynomial in $\fq[x]$ with degree $2$ and order $2^{L_2}5$.  

Moreover, because $2 \nmid (q^2-1)/2^{L_2}5$, by Lemma~\ref{irred},  $x^4 - w \rho_{L_1}$ is also irreducible. Hence 
\[
Q_{2^{L_4} 5}(x) = Q_{2^{L_2} 5}(x^2) = \prod_{w \in \Omega(5)} \prod_{\rho_{L_1} \in \Omega(2^{L_1})} \left( x^4 - w \rho_{L_1} \right).
\]
In this case, $x^4 - w  \rho_{L_1}$ is an irreducible polynomial of degree $4$ and
order $2^{L_4} 5$. In general, for $n > L_4$, then $2^{n-L_1} \equiv 0 \pmod{4}$ and $q^{4} -1 \equiv 1 \pmod{4}$. Because $2^{n-L_1} \nmid (q^4 -1)/2^{L_4}5$, by Lemma~\ref{irred},  the polynomials $x^{2^{n-L_1}} - w  \rho_{L_1}$ are irreducible over $\fq$  and thus
\[
Q_{2^n 5}(x) = Q_{2^{L_1} 5}(x^{2^{n-L_1}}) = \prod_{w \in \Omega(5)} \prod_{\rho_{L_1} \in \Omega(2^{L_1})}   \left( x^{2^{n-L_1}} - w  \rho_{L_1} \right). 
\]
\rmv{
If $k$ is odd,   then $L_1 = v_2(q-1) =2$ and thus $\rho$ is a non-square element in $\fq$. Hence $\rho\rho_{L_1}$   is a square element (either $\rho_2^2$ or $-\rho_2^2$) in $\fq$. Then $x^2 + w \rho \rho_{L_1}$ is not irreducible polynomial in $\fq[x]$. Hence
\begin{eqnarray*}
Q_{2^{L_2} 5}(x) &=& Q_{2^{L_1} 5}(x^2) =  Q_{2^2 5}(x^2) \\
&=& \prod_{w \in \Omega(5)} \prod_{\rho_{L_1} \in \Omega(2^{L_1})}   \left( x^2  +  w \rho \rho_{L_1} \right) \\
&=& \prod_{w \in \Omega(5)}  \left( x^2  +  w \rho_2^2 \right)  \left( x^2  -  w \rho_2^2 \right) \\
& =& \prod_{u \in \Omega(5)}   \left( x  +  u \rho \rho_{2} \right)  \left( x  -  u \rho \rho_{2} \right)  \left( x  +  u \rho_{2} \right)  \left( x  - u \rho_{2} \right)\\
& =& \prod_{u \in \Omega(5)}   \left( x  +  u \rho_{1} \right)  \left( x  -  u  \rho_{1} \right)  \left( x  +  u \rho_{2} \right)  \left( x  - u \rho_{2} \right)
\end{eqnarray*}

 Similarly, we obtain, for $n > L_2$, that $x^{2^{n-L_2}} + w \rho \rho_{L_1}$ is irreducible
\[
Q_{2^n 5}(x) = Q_{2^{L_2} 5}(x^{2^{n-L_2}}) = \prod_{w \in \Omega(5)} \prod_{\rho_{L_1} \in \Omega(2^{L_1})}   \left( x^{2^{n-L_1-1}} + w  \rho \rho_{L_1} \right)\left( x^{2^{n-L_1-1}} - w  \rho \rho_{L_1} \right). 
\]
}
\end{proof}

\begin{theorem}
Let  $q = 20k + 11$  for some nonnegative integer $k$. 
Then we have the following factorization of $2^n5$-th cyclotomic polynomial $Q_{2^n 5} (x)$  over $\fq$.

(i) For $n=0, 1, 2$, we have 
\[
Q_{5} (x) = \prod_{w \in \Omega(5)}   \left( x - w  \right),  
~
Q_{10} (x) = \prod_{w \in \Omega(5)}   \left( x + w\right)
,
~
Q_{20} (x) = \prod_{w \in \Omega(5)}  \left( x^2 + w  \right).
\]

(ii)  if $n \geq L_2 = 3$, then 
\[
Q_{2^n 5} (x) =
\left\{ 
\begin{array}{rr} 
 \displaystyle{\prod_{w \in \Omega(5)}}  \prod_{c^2 = -2} \left( x^{2^{n-2}} + cw x^{2^{n-3}}  -w^2 \right) , &
if~ $k$ ~is ~even;  \\
& \\
\displaystyle{\prod_{w \in \Omega(5)}}  \prod_{c^2 = 2} \left( x^{2^{n-2}} + cw x^{2^{n-3}}  + w^2 \right),  &
if~$k$ ~is~ odd; 
\end{array}
\right.
\]

\end{theorem}

\begin{proof} In this case, $L_1 = v_2(q-1) = 1$, $L_2 = v_2(q^2-1) =3$, and
$L_4 = L_2 +1 = 4$. 
It is obvious that $Q_{5} (x) = \prod_{w \in \Omega(5)}   \left( x - w  \right)$, and $Q_{10} (x) =  Q_5(-x) = \prod_{w \in \Omega(5)}   \left( x + w\right)$ because $5 \mid q-1$. Because $q\equiv 3 \pmod{4}$, $-1$ is a non-square in $\fq$.  Hence $x^2 + w$ is irreducible in $\fq[x]$ and then $Q_{20} (x) = Q_{10}(x^2) = \prod_{w \in \Omega(5)}   \left( x^2 + w\right)$.  
Now $Q_{40}(x) = Q_{20}(x^2) =  \prod_{w \in \Omega(5)}   \left( x^4 + w\right)$.  Because $\gcd(40, q) =1$ and $q^2 \equiv 1 \pmod{40}$, by Theorem 2.47 in \cite{LN}, $Q_{40}(x)$ factors into $\phi(40)/2 = 12$ distinct monic quadratic irreducible polynomials in $\fq[x]$. Hence $x^4 + w$ can be factorized into a product of two monic quadratic polynomials. Let $x^4 + w = (x^2 + ax + b)(x^2 +cx+d)$ where $a, b, c, d\in \fq$. Comparing both sides,  we have $a+c =0$, $b+d+ac=0$, $bd+ad=0$, and $bd =w$. Replacing $a$ by $-c$, we obtain $b+d -c^2 =0$, $(b-d)c =0$, and $bd=w$.  Because $w^5 =1$, we can write $w=v^4$ where $v = w^{-1}$. Because $-1$ is a non-square,  we can only have two possible solutions ($c\neq 0$ because, otherwise, $b^2 = -v^4$, a contradiction):
(i) $b=d =v^2$, $a=-c$, and $c^2 = 2v^2$ if $2$ is a square; (ii) $b=d =-v^2$, $a=-c$, and $c^2=-2v^2$ if $-2$ is a square.  We note that $q \equiv 4k+3 \pmod{8}$. Hence if $k$ is even, then $q\equiv 3 \pmod{8}$; otherwise, $q\equiv 7 \pmod{8}$.Moreover, $q\equiv 3 \pmod{8}$ implies that the characteristic $p$ of $\fq$ also satisfies $p \equiv 3 \pmod{8}$, therefore $-2$ is a square in $\fq$ if $k$ is even. Similarly, $2$ is a square in $\fq$ if $k$ is odd.  As $w$ ranges over $\Omega(5)$, $v$ also ranges over $\Omega(5)$. Hence we obtain
\[
Q_{40} (x) =
\left\{ 
\begin{array}{rr} 
 \displaystyle{\prod_{w \in \Omega(5)}}  \prod_{c^2 = -2} \left( x^{2} + cw x  -w^2 \right) , &
if~ $k$ ~is ~even;  \\
& \\
\displaystyle{\prod_{w \in \Omega(5)}}  \prod_{c^2 = 2} \left( x^{2} + cw x  + w^2 \right),  &
if~$k$ ~is~ odd; 
\end{array}
\right.
\]
Each $x^2 +cwx-w^2$ is an irreducible polynomial of degree $2$ and order $40$.
If $n \geq 4$, then $2^{n-2} \equiv 0 \pmod{4}$ and $q^2 -1 \equiv 0 \pmod{4}$. Moreover, $2^{n-2} \nmid (q^2 -1)/40$. Hence by Lemma~\ref{irred},  we have that
$x^{2^{n-2}} +cw x^{2^{n-3}} -w^2$ is irreducible and 
\[
Q_{2^n 5} (x) = Q_{2^3 5}(x^{2^{n-3}}) =
\left\{ 
\begin{array}{rr} 
 \displaystyle{\prod_{w \in \Omega(5)}}  \prod_{c^2 = -2} \left( x^{2^{n-2}} + cw x^{2^{n-3}}  -w^2 \right) , &
if~ $k$ ~is ~even;  \\
& \\
\displaystyle{\prod_{w \in \Omega(5)}}  \prod_{c^2 = 2} \left( x^{2^{n-2}} + cw x^{2^{n-3}}  + w^2 \right),  &
if~$k$ ~is~ odd; 
\end{array}
\right.
\]
\end{proof}

\begin{theorem}
Let $q\equiv  9 \pmod{20}$ and $L_i = v_2(q^i-1)$ for $i \geq 1$. 
Then we have the following factorization of $2^n 5$-th cyclotomic polynomial $Q_{2^n 5} (x)$  over $\fq$.

(i) For $n=0, 1$, we have 
\[
Q_{5} (x) = \prod_{a = w + w^{-1} \atop{w \in \Omega(5)}}   \left( x^2 - a x + 1  \right),  
~
Q_{10} (x) = \prod_{a = w + w^{-1} \atop{w \in \Omega(5)}}   \left( x^2 + a x + 1  \right).
\] 

(ii) if $2 \leq n < L_2$, then 
\[
Q_{2^n 5} (x) = \prod_{a_n = \rho_{n} (w + w^{-1}) \atop{w \in \Omega(5)}} \prod_{\rho_{n} \in \Omega(2^{n})}  \left( x^2 + a_n x +  \rho_{n-1} \right).
\]

(iii) if $n \geq L_2 $, then 
\[
Q_{2^n 5} (x) = \prod_{a_{L_2}^2 = 2\rho_{L_1} - a_{L_1}} \prod_{ a_{L_1} = \rho_{L_1}(w + w^{-1}) \atop{w \in \Omega(5)}} \prod_{\rho_{L_1} \in \Omega(2^{L_1})}  \left( x^2 + a_{L_2} x +  \rho_{L_1} \right)
\]
\end{theorem}

\begin{proof}
Let $q =20k + 9$. If $k$ is odd, then $L_1 =2$. If $k$ is even, then $L_1 =3$. 
In this case,  $L_2 = L_1+1$, and $L_4 = L_2 +1$. 
Because $5 \nmid q-1$,  $w \in \Omega(5)$ implies $w \not\in \fq$. However, $5 \mid q+1$ implies $a=w+w^{-1} = w+w^{q}  \in \fq$. Hence 
$ Q_{5}(x) =  \displaystyle{ \prod_{a = w + w^{-1} \atop{w \in \Omega(5)}}   \left( x^2 - a x + 1  \right) }$ and 
$Q_{10} (x) = Q_5(-x) =  \displaystyle{ \prod_{a = w + w^{-1} \atop{w \in \Omega(5)}}   \left( x^2 + a x + 1  \right)} = \displaystyle{ \prod_{a = w + w^{-1} \atop{w \in \Omega(5)}}   \left( x^2 + a x + \rho_0 \right)} $.

Then
$Q_{20}(x) = Q_{10}(x^2) = \displaystyle{ \prod_{a = w + w^{-1} \atop{w \in \Omega(5)}}   \left( x^4 + a x^2 + \rho_{0}  \right)}$. 
Again, $q^2 \equiv 1 \pmod{2^{L_1} 5}$ and Theorem~2.47 in \cite{LN} imply that $Q_{2^n 5}(x)$ factors into distinct monic quadratic polynomials. Let $a_1 =a$.  In order to factor $Q_{20}(x)$,  we need to factor $x^4 + a_{1} x^2 + \rho_{0}$ into monic quadratic irreducible polynomials. 
Let $x^4 + a_{1}x^2 + \rho_{0} = (x^2 + b x+c)(x^2 + dx + e)$ where $b, c, d, e \in \fq$. Then we obtain
\begin{eqnarray*}
b+d &=&0\\
c+e+bd &=& a_{1}\\
be+cd &=&0\\
ce &=& \rho_{0}
\end{eqnarray*}
Hence $b=-d$. Continue to solve the above system,  either we have $b=-d =0$ or $e=c$. If $b=-d =0$, then $c$ satisfies that $c^2 -ac+1 =0$, contradicts to that $x^2 -ax +1$ is irreducible. Hence $e=c$. Let $x^4 + a_{1} x^2 +\rho_{0} = (x^2 + a_{2} x + c)(x^2-a_2 x+ c)$. 
Therefore $e=c = \pm \rho_{1} \in \fq$ and $a_2^2 = \pm 2\rho_{1} -a_{1}$. Since we can verify directly that 
$a_2=  \rho_{2} (w^3+w^{-3}) \in \fq$ are solutions to $a_2^2 =  2\rho_{1} -a_{1}$ where $\rho_2^2 = \rho_1$, we obtain
\[
Q_{20} (x) = \prod_{a_2 =  \rho_{2} (w + w^{-1}) \atop{w \in \Omega(5)}} \prod_{\rho_{2} \in \Omega(2^{2})}  \left( x^2 +  a_2 x +  \rho_{1} \right).
\]

If $k$ is odd, then
\[
Q_{40} (x) = Q_{20}(x^2) = \prod_{a_2 =  \rho_{2} (w + w^{-1}) \atop{w \in \Omega(5)}} \prod_{\rho_{2} \in \Omega(2^{2})}  \left( x^4 +  a_2 x^2 +  \rho_{1} \right).
\]
Similarly, let $\left( x^4 +  a_2 x^2 +  \rho_{1} \right)=\left( x^2 +  a_{3} x +  \rho_{2} \right) 
\left( x^2 -  a_3 x +  \rho_{2} \right)$ where $a_3, \rho_2 \in \fq$. Hence $a_3^2 = 2\rho_2 - a_2$. Because
$ 2\rho_2 -\rho_2(w+w^{-1})  = - \rho_2 (w^3-w^{-3})^2 \in \fq$ and both $\rho_2$ and $- (w^3-w^{-3})^2$ are non-square elements in $\fq$,  there exist $a_3 \in \fq$ such that $a_3^2 = 2\rho_2 -a_2$. Hence
\[
Q_{40} (x) = 
\prod_{a_3^2 = 2 \rho_2 - a_2} \prod_{a_2 =  \rho_{2} (w + w^{-1}) \atop{w \in \Omega(5)}} \prod_{\rho_{2} \in \Omega(2^{2})}  \left( x^2 +  a_3 x +  \rho_{2} \right).
\]
Moreover, for any $n\geq 4$, by Theorem~\ref{general}, we conclude $x^{2^{n-2}} + a_3 x^{2^{n-3}} + \rho_2$ is irreducible. Hence 
\[
Q_{2^n 5} (x) = 
\prod_{a_3^2 = 2 \rho_2 - a_2} \prod_{a_2 =  \rho_{2} (w + w^{-1}) \atop{w \in \Omega(5)}} \prod_{\rho_{2} \in \Omega(2^{2})}  \left( x^{2^{n-2}} +  a_3 x^{2^{n-3}} +  \rho_{2} \right).
\]

If $k$ is even, then $L_1 =3$ and $\rho_3 \in \fq$. Then 
\[
Q_{40} (x) = Q_{20}(x^2) = \prod_{a_3 =  \rho_{3} (w + w^{-1}) \atop{w \in \Omega(5)}} \prod_{\rho_{3} \in \Omega(2^{3})}  \left( x^2 +  a_3 x +  \rho_{3} \right).
\]
Similarly, fr $n \geq 4$, we have
\[
Q_{2^n 5} (x) = 
\prod_{a_4^2 = 2 \rho_3 - a_3} \prod_{a_3 =  \rho_{3} (w + w^{-1}) \atop{w \in \Omega(5)}} \prod_{\rho_{3} \in \Omega(2^{3})}  \left( x^{2^{n-2}} +  a_4 x^{2^{n-3}} +  \rho_{3} \right).
\]

\end{proof}

\begin{theorem}
Let $q\equiv  19 \pmod{20}$. 
Then we have the following factorization of $2^n 5$-th cyclotomic polynomial $Q_{2^n 5} (x)$  over $\fq$.

(i) For $n=0, 1$, we have 
\[
Q_{5} (x) = \prod_{a = w + w^{-1} \atop{w \in \Omega(5)}}   \left( x^2 - a x + 1  \right),  ~ 
Q_{10} (x) = \prod_{a = w + w^{-1} \atop{w \in \Omega(5)}}   \left( x^2 + a x + 1  \right).
\] 

(ii) 
\[
Q_{20} (x) =\prod_{\rho_{2} \in \Omega(2^{2})} \prod_{a_2 = \rho_2 (w - w^{-1}) \atop{w \in \Omega(5)}}   \left( x^2 + a_2 x +  1 \right).
\]

(iii)  if $n \geq  L_2 = 3$, then 
\[
Q_{2^n 5} (x) =\prod_{\rho_{3} \in \Omega(2^{3})}  \prod_{a_3 = \rho_2 \rho_3w - (\rho_2 \rho_3w)^{-1} \atop{w \in \Omega(5)}}  \left( x^{2^{n-2}} + a_3 x^{n-3} -  1 \right).
\]

\end{theorem}

\begin{proof}
In this case, $L_1 = 1$, $L_2 = L_1+2 = 3$, and $L_4 = L_2 + 1$. 
Again, $5\nmid q-1$ implies that if $w\in \Omega(5)$ then $w\not\in \fq$. However, $w\in \F_{q^2}$. Moreover, $-1$ is a non-square element. 
Again, it is trivial to obtain 
\[
Q_{5} (x) = \prod_{a = w + w^{-1} \atop{w \in \Omega(5)}}   \left( x^2 - a x + 1  \right),  ~ 
Q_{10} (x) = \prod_{a = w + w^{-1} \atop{w \in \Omega(5)}}   \left( x^2 + a x + 1  \right).
\] 
Let $a_1$ denotes $w+w^{-1}$. 
Because $q\equiv 19 \pmod{20}$,   we have $(\rho_2w)^{q+1} = 1$ and thus $\rho_2 w + (\rho_2 w)^{-1}  = \rho_2 w + (\rho_2 w)^q \in \fq$. Moreover,
\[
x^4 + a_1 x^2 + 1 = (x^2 + a_2 x + 1)(x^2-a_2x+1),
\]
where $a_2 = \rho_2 (w^3-w^{-3}) = \rho_2 w^3 + (\rho_2 w)^{-3} \in \fq$. Again $w$ ranges over $\Omega(5)$ means that $w^3$ ranges over $\Omega(5)$. Therefore,
\[
Q_{20} (x) =\prod_{\rho_{2} \in \Omega(2^{2})} \prod_{a_2 = \rho_2 (w - w^{-1}) \atop{w \in \Omega(5)}}   \left( x^2 + a_2 x +  1 \right).
\]

For $n=3$, let $\rho_3^2 = \rho_2$ and  $a_3 = \rho_2 \rho_3 w^3 - (\rho_2 \rho_3 w^3)^{-1}$. We claim that
that  $a_3 \in \fq$. First, we note that $\rho_2^q = \rho_2^{-1}$,  and $w^q = w^{-1}$. Moreover, $\rho_3^{q}  = - \rho_3^{-1}$ because $\rho_3^{2(q+1)} =1$
and $\rho_3^{q+1} \neq 1$. Then
\[
a_3^q = (\rho_2 \rho_3 w^3 - (\rho_2 \rho_3 w^3)^{-1})^q = \rho_2^q \rho_3^q w^{3q} - \rho_2^{-q} \rho_3^{-q} w^{-3q} =  (\rho_2^{-1}) (-\rho_3^{-1}) w^{-3} - (\rho_2 (-\rho_3) w^3) = a_3
\]
Moreover, $ (x^2 +a_3 x -1)(x^2-a_3 x-1) = x^4 + a_2 x + 1$ because $a_3^2 = -\rho_2 w + \rho_2^{-1} w^{-1} -2 =  -a_2 -2$. Hence
\[
Q_{40} (x) =\prod_{\rho_{3} \in \Omega(2^{3})}  \prod_{a_3 = \rho_2 \rho_3w - (\rho_2 \rho_3w)^{-1} \atop{w \in \Omega(5)}}  \left( x^{2} + a_3 x -  1 \right).
\]
The rest of proof follows from Theorem~\ref{general} and $Q_{2^n 5}(x) = Q_{40}(x^{2^{n-3}})$. 
\end{proof}

\section{Case: $q\equiv \pm 2 \pmod{5}$ and $q\equiv 1 \pmod{4}$}
\label{new1}

We note that if $q\equiv \pm 2 \pmod 5$ and $q\equiv 1 \pmod 4$ (i.e., $q \equiv 13 \pmod{20}$ or $q\equiv 17 \pmod{20}$), then $L_4 = L_2+1$ and $L_2 = L_1 + 1$. Moreover, $\rho_1 = -1$ must be a square and thus there exists $\rho_2\in \fq$ such that $\rho_2^2 = \rho_1$.

\begin{theorem} \label{factor-mod-13-17}
Let $q\equiv \pm 2 \pmod 5$ and $q\equiv 1 \pmod 4$. 
Then we have the following factorization of $2^n 5$-th cyclotomic polynomial $Q_{2^n 5} (x)$  over $\fq$.

(i) If $0\leq n \leq L_1$, then 
\[
Q_{2^n 5} (x) = \prod_{\rho_n \in \Omega(2^n)}  \left( x^4 + \rho_n x^3 + \rho_n^2 x^2 + \rho_n^3 x + \rho_n^4\right).
\]

(ii) If $n = L_2$ (i.e., $L_2 = L_1 + 1$), then 
\[
Q_{2^n 5} (x) = \prod_{\rho_{n-1} \in \Omega(2^{n-1})} \prod_{a_{n}^2 = 5 \rho_{n-1}}  \left( x^4 + a_n x^3 + 3\rho_{n-1} x^2 + a_n \rho_{n-1} x + \rho_{n-2}\right).
\]

(iii) If $ n = L_4$ (i.e., $L_4 = L_2 + 1$), then 
\[
Q_{2^n 5} (x) = \left\{ 
\begin{array}{rr} 
 \displaystyle{\prod_{\rho_{n-2} \in \Omega(2^{n-2})} \prod_{a_{n-1}^2 = 5 \rho_{n-2}} \prod_{a_n^2 = (2\rho_2-1)a_{n-1}} } \left( x^4 + a_n x^3 + a_{n-1}\rho_{2} x^2 + (-5\rho_{n-2})a_n^{-1}  x - \rho_{n-2}\right), \\
if~ (2\rho_{2} -1)a_{n-1} ~is~ a~square; & \\
& \\
\displaystyle{\prod_{\rho_{n-2} \in \Omega(2^{n-2})} \prod_{a_{n-1}^2 = 5 \rho_{n-2}} \prod_{a_n^2 = -(2\rho_2+1)a_{n-1}} } \left( x^4 + a_n x^3 + (-a_{n-1}\rho_{2}) x^2 + (-5\rho_{n-2})a_n^{-1}  x - \rho_{n-2}\right), \\
 if~ (2\rho_{2} -1)a_{n-1} ~is~ a~nonsquare. & 
\end{array}
\right.
\]

(iv) If $n > L_4=L$, then $Q_{2^n 5} (x)$ can be factorized as 
\[
 \left\{ 
\begin{array}{rr} 
 \displaystyle{\prod_{\rho_{L_1} \in \Omega(2^{L_1})} \prod_{a_{L_2}^2 = 5 \rho_{L_1}} \prod_{a_{L_4}^2 = (2\rho_2-1)a_{L_2}} } \left( x^{2^{n-L_4+2}} + a_{L_4} x^{3\cdot 2^{n-L_4}} + a_{L_2}\rho_{2} x^{2^{n-L_4+1}} + (-5\rho_{L_1})a_{L_4}^{-1}  x^{2^{n-L_4}} - \rho_{L_1}\right), \\
if~ (2\rho_{2} -1)a_{L_2} ~is~ a~square; & \\
& \\
\displaystyle{\prod_{\rho_{L_1} \in \Omega(2^{L_1})} \prod_{a_{L_2}^2 = 5 \rho_{L_1}} \prod_{a_{L_4}^2 = -(2\rho_2+1)a_{L_2}} } 
\left( x^{2^{n-L_4+2}} + a_{L_4} x^{3\cdot 2^{n-L_4}} + (-a_{L_2}\rho_{2}) x^{2^{n-L_4+1}} + (-5\rho_{L_1})a_{L_4}^{-1}  x^{2^{n-L_4}} - \rho_{L_1}\right), \\
 if~ (2\rho_{2} -1)a_{L_2} ~is~ a~ nonsquare. & 
\end{array}
\right.
\]

\end{theorem}

\begin{proof}
Because the smallest positive $d$ satisfying $q^d \equiv 1 \pmod{5}$ is $4$ under our assumption, by Theorem 2.47 in \cite{LN}, for all $ 0 \leq n \leq L_4$,  $Q_{2^n 5}$ factors into a product of $\phi(2^n 5)/ 4 = 2^{n-1}$ distinct monic irreducible polynomials of degree $4$.

(i) If $n \leq L_1$, then $\rho_n \in \Omega (2^n) \subseteq \fq$. Hence $x^4 + \rho_n x^3 + \rho_n^2 x^2 + \rho_n^3 x + \rho_n^4 \in \fq[x]$. The factorization of $Q_{2^n 5}(x)$ when $n=0, 1$ is trivial. Moreover, it is straightforward to verify that for $\rho_n^2 = \rho_{n-1}$
\[
x^8 + \rho_{n-1} x^6 + \rho_{n-1}^2 x^4 + \rho_{n-1}^3 x^2 + \rho_{n-1}^4
= \left( x^4 + \rho_n x^3 + \rho_n^2 x^2 + \rho_n^3 x + \rho_n^4 \right)\left(
x^4 - \rho_n x^3 + \rho_n^2 x^2 - \rho_n^3 x + \rho_n^4 \right).
\]
Hence (i) follows from  $Q_{2^n 5}(x) = Q_{2^{n-1} 5}(x^2)$ and the consequences of Theorem 2.47 in \cite{LN} as mentioned above.

(ii) From (i) and Lemma~\ref{cyclo},  we have  
\[
Q_{2^{L_2} 5}(x) = Q_{2^{L_1} 5}(x^2) =
\prod_{\rho_{L_1} \in \Omega(2^{L_1})}  \left( x^8 + \rho_{L_1} x^6 + \rho_{L_1}^2 x^4 + \rho_{L_1}^3 x^2 + \rho_{L_1}^4\right).
\]
Because the Legendre symbol $ \left( 5 \atop p \right) = 1$ iff $p \equiv \pm 1, \pm 9 \pmod{20}$, $5$ is a non-square in $\fq$ under the assumption of our theorem. Hence $5 \rho_{L_1}$ is a square element in $\fq$ as $\rho_{L_1}$ is also a non-square element in $\fq$. Let $a_{L_2}^2 = 5 \rho_{L_1}$. Then $\pm a_{L_2} \in \fq$. Hence $x^4 \pm  a_{L_2} x^3 + 3 \rho_{L_1} x^2 \pm a_{L_2} \rho_{L_1} x + \rho_{L_1-1} \in \fq[x]$. One can also easily verify that
\[
x^8 + \rho_{L_1} x^6 + \rho_{L_1}^2 x^4 + \rho_{L_1}^3 x^2 + \rho_{L_1}^4
= \prod_{a_{L_2}^2 = 5\rho_{L_1}} \left( x^4 + a_{L_2} x^3 + 3 \rho_{L_1} x^2 + a_{L_2} \rho_{L_1} x + \rho_{L_1-1} 
 \right).
\]
Therefore the rest of proof of (ii) follows.

(iii) In this case, we essentially need to factor 
$
 x^8 + a_{L_2} x^6 + 3 \rho_{L_1} x^4 + a_{L_2} \rho_{L_1} x^2 + \rho_{L_1-1} 
$
into two monic quartic polynomials in $\fq[x]$, where $\rho_{L_1} \in \Omega(2^{L_1})$ and $a_{L_2}^2 = 5\rho_{L_1}$. 

Because $-1$ is a square element and $5$ is a non-square, $-(2\rho_2 + 1)a_{L_2} (2\rho_2 -1) a_{L_2} = -(4\rho_2^2 -1) a_{L_2}^2 = 5 a_{L_2}^2$ is a non-square element in $\fq$.  Hence either $-(2\rho_2 + 1)a_{L_2}$ or $(2\rho_2 -1) a_{L_2}$ (exactly one of them) is a square element in $\fq$. 

If $(2\rho_2 -1) a_{L_2}$ is a square, we let $a_{L_4}^2 = (2\rho_2 -1) a_{L_2}$. Then 
\[
 x^8 + a_{L_2} x^6 + 3 \rho_{L_1} x^4 + a_{L_2} \rho_{L_1} x^2 + \rho_{L_1-1} 
=
\prod_{a_{L_4}^2 = (2\rho_2 -1) a_{L_2}} \left(
x^4 + a_{L_4} x^3 + \rho_2 a_{L_2} x^2 + (-5\rho_{L_1}) a_{L_4}^{-1} x - \rho_{L_1}
\right)
\]

If $(2\rho_{2} -1)a_{n-1}$ is a non-square then $-(2\rho_2 + 1)a_{L_2}$ is a square. In this case, we let $a_{L_4}^2 = -(2\rho_2 + 1) a_{L_2}$. Then
\[
 x^8 + a_{L_2} x^6 + 3 \rho_{L_1} x^4 + a_{L_2} \rho_{L_1} x^2 + \rho_{L_1-1} 
=
\prod_{a_{L_4}^2 = -(2\rho_2 + 1) a_{L_2}} \left(
x^4 + a_{L_4} x^3 - \rho_2 a_{L_2} x^2 + (-5\rho_{L_1}) a_{L_4}^{-1} x - \rho_{L_1}
\right)
\]
Because each quartic polynomial is in $\fq[x]$ and $Q_{2^{L_4} 5}(x)$ factors into product of quartic polynomials, every such quartic polynomial must be irreducible. Hence (iii) is proved.

(iv) For $n \geq L_4$, then $Q_{2^n 5}(x) = Q_{2^{L_4} 5} (x^{2^{n-L_4}})$. 
If $(2\rho_2 -1) a_{L_2}$ is a square, then  each irreducible factor 
\[
x^4 + a_{L_4} x^3 + \rho_2 a_{L_2} x^2 + (-5\rho_{L_1}) a_{L_4}^{-1} x - \rho_{L_1}
\] of $Q_{2^{L_4} 5}(x)$ has degree $4$ and order $2^{L_4}5$, where $a_{L_4}^2 = (2\rho_2 -1) a_{L_2}$.  By Theorem~\ref{general},  
\[
 x^{2^{n-L_4+2}} + a_{L_4} x^{3\cdot 2^{n-L_4}} + a_{L_2}\rho_{2} x^{2^{n-L_4+1}} + (-5\rho_{L_1})a_{L_4}^{-1}  x^{2^{n-L_4}} - \rho_{L_1}
\] must also irreducible.

 Similarly, if $(2\rho_2 -1) a_{L_2}$ is a non-square, then
\[
 x^{2^{n-L_4+2}} + a_{L_4} x^{3\cdot 2^{n-L_4}} - a_{L_2}\rho_{2} x^{2^{n-L_4+1}} + (-5\rho_{L_1})a_{L_4}^{-1}  x^{2^{n-L_4}} - \rho_{L_1}
\] where $a_{L_4}^2 = -(2\rho_2 + 1) a_{L_2}$ must be irreducible. Hence the proof is complete.
\end{proof}

\begin{corol} \label{irred-mod-13-17}
Let $q\equiv \pm 2 \pmod 5$ and $q\equiv 1 \pmod 4$. 
Let $\rho_n \in \Omega(2^n)$, $\rho_n^{2^{n-2}} = \rho_2$, and $a_{L_2}^2 = 5 \rho_{L_1}$.

(i) If $2(\rho_{2} -1)a_{L_2}$ is a square, 
 then $x^{2^{n-L_4+2}} + a_{L_4} x^{3\cdot 2^{n-L_4}} + a_{L_2}\rho_{2} x^{2^{n-L_4+1}} + (-5\rho_{L_1})a_{L_4}^{-1}  x^{2^{n-L_4}} - \rho_{L_1}$ is irreducible over $\fq$ for each choice of $\rho_n$, $a_{L_2}$,  and $a_{L_4}^2 = (2\rho_2 -1)a_{L_2}$;

(ii)  otherwise, 
$x^{2^{n-L_4+2}} + a_{L_4} x^{3\cdot 2^{n-L_4}} + (-a_{L_2}\rho_{2}) x^{2^{n-L_4+1}} + (-5\rho_{L_1})a_{L_4}^{-1}  x^{2^{n-L_4}} - \rho_{L_1}$ is irreducible over $\fq$ for each choice of  $\rho_n$, $a_{L_2}$,  and  $a_{L_4}^2 = -(2\rho_2 +1)a_{L_2}$.  
\end{corol}

\section{Case: $q\equiv \pm 2 \pmod{5}$ and $q\equiv 3 \pmod{4}$ }
\label{new2}

First we note that both $-1$ and $5$ are non-square elements and thus $-5$ is a square element in $\fq$. 
We note also that if $q\equiv \pm 2 \pmod 5$ and $q\equiv 3 \pmod 4$ (i.e., $q \equiv 3 \pmod{20}$ or $q\equiv 7 \pmod{20}$), then $L_4 = L_2+1$ and $L_1 =1$. However, $L_2 = L_1 + 2$ for $q\equiv 3 \pmod{20}$ and $L_2 = L_1 + 2$ or $L_1+3$ for $q\equiv 7 \pmod{20}$, depending on parity of $k$ when $q=20k+7$. 
Hence we separate this case into two subcases.

\subsection{Case $q \equiv 3 \pmod{20}$}
\label{new2}

In this case, $L_1 = 1$, $L_2 = L_1 + 2 =3$ and $L_4 = L_2+1 =4$.

\begin{theorem} \label{factor-mod-3}
Let $q\equiv 3 \pmod{20}$ and $q > 3$.  Then we have the following factorization of cyclotomic polynomials $Q_{2^n 5} (x)$  
over $\fq$.

(i) If $0 \leq n \leq 1$, then  
\[
Q_{2^n 5} (x) = \prod_{\rho_n \in \Omega(2^n)}  \left( x^4 + \rho_n x^3 + \rho_n^2 x^2 + \rho_n^3 x + \rho_n^4\right).
\]

(ii) If $n = 2$, then 
\[
Q_{2^n 5} (x) =  \prod_{a_{2}^2 = -5 }  \left( x^4 + a_2 x^3 + 3\rho_{1} x^2 + a_2 \rho_1 x + 1\right).
\]

(iii) If $n = 3$, then 
\[
Q_{2^n 5} (x) =
 \left\{ 
\begin{array}{ll} 
 \displaystyle{ \prod_{a_{2}^2 = -5} \prod_{ a_3^2 = 2 b_3 - a_{2} \atop{b_3= 1} } } \left( x^4 + a_3 x^3 + b_3 x^2 + 3a_3^{-1}  x + 1 \right), &
if~ 2  - a_{2} ~is~ a~ square;  \\
& \\
\displaystyle{ \prod_{a_{2}^2 = -5} \prod_{a_3^2 = 2 b_3 - a_{2} \atop{b_3= - 1} } } \left( x^4 + a_3 x^3 + b_3 x^2 + 3a_3^{-1}  x + 1 \right), &
 if~ -2 - a_{2} ~is~a~ square. 
\end{array}
\right.
\]

(iv) If $ n = 4$, then 

\[
Q_{2^n 5} (x) = \left\{ 
\begin{array}{ll} 
 \displaystyle{ \prod_{a_{2}^2 = -5} \prod_{ a_3^2 = 2  - a_{2} \atop{c_3 = 3a_3^{-1}}} \prod_{a_4, b_4, c_4} } \left( x^4 + a_4 x^3 + b_4 x^2 + c_4  x - 1 \right), & 
if~ 2 - a_{2} ~is~ a~square; \\
& \\
\displaystyle{ \prod_{a_{2}^2 = -5} \prod_{a_3^2 = - 2  - a_{2}\atop{c_3 = 3a_3^{-1}} }  
\prod_{a_4, b_4, c_4} }
\left( x^4 + a_4 x^3 + b_4 x^2 + c_4 x + 1 \right), &
 if~ - 2 - a_{2} ~is~ a~ square,
\end{array}
\right.
\]
where $a_4, b_4, c_4$ satisfy either 
\begin{equation}\label{cond4a}
a_4^2 =2b_4 -a_3, b_4 = \alpha + a_2~ or ~\alpha + a_2, \alpha^2 = -2, c_4^2 = -2b_4-c_3, ~and~ c_4 = (b_4^2 -3)(2a_4)^{-1}, 
\end{equation} 
when $2-a_2$ is a square, or
\begin{equation} \label{cond4b}
a_4^2 =2b_4 -a_3,\\ b_4 =  \beta +1 ~or~\beta-1 , \beta^2 = 2, \\ c_4^2 = 2b_4-c_3, ~and~ c_4 = (b_4^2 +3)(2a_4)^{-1}, 
\end{equation}
when $-2-a_2$ is a square.

(v) If $n > 4$, then $Q_{2^n 5} (x)$ can be factorized as 

\[
\left\{ 
\begin{array}{ll} 
 \displaystyle{ \prod_{a_{2}^2 = -5} \prod_{ a_3^2 = 2  - a_{2} \atop{c_3 = 3a_3^{-1}} } \prod_{a_4, b_4, c_4 } } \left( x^{2^{n-2}} + a_4 x^{3\cdot 2^{n-4}} + b_4 x^{2^{n-3}} + c_4 x^{2^{n-4}} - 1 \right), & 
if~ 2 - a_{2} ~is~ a~square; \\
& \\
\displaystyle{ \prod_{a_{2}^2 = -5} \prod_{a_3^2 = - 2  - a_{2} \atop{c_3 = 3a_3^{-1}}}  
\prod_{ a_4, b_4, c_4} }
\left( x^{2^{n-2}} + a_4 x^{3\cdot 2^{n-4}} + b_4 x^{2^{n-3}} + c_4  x^{2^{n-4}} + 1 \right), &
 if~ - 2 - a_{2} ~is~ a~ square,  
\end{array}
\right.
\]
where $a_4, b_4, c_4$ satisfy either (\ref{cond4a}) if $2-a_2$ is a square,  or (\ref{cond4b}) if $-2-a_2$ is a square. 
\end{theorem}

\begin{proof}

(i) Again $Q_{2^n 5}$ factors into a product of $\phi(2^n 5)/ 4 = 2^{n-1}$ distinct monic irreducible polynomials of degree $4$.  It is trivial to check $Q_5(x) = \prod_{w\in \Omega(5)} (x-w) = x^4 - x^3 + x^2 - x + 1$ and $Q_{10}(x) = Q_{5}(-x)= x^4 + x^3 + x^2 + x +1$.

(ii) Because $-5$ is a square when $q\equiv 3 \pmod{20}$  and $q >3$, for each
$a_2$ satisfying $a_2^2 = -5$, we can verify that
\[
\left( x^4 + a_2 x^3 + 3 \rho_1 x^2 + a_2 \rho_1 x +1 \right)
\left( x^4 - a_2 x^3 + 3 \rho_1 x^2 - a_2 \rho_1 x +1 \right)
= x^8 + x^6 + x^4 + x^2 + 1 
\]
Hence the proof follows from $Q_{20}(x) = Q_{10}(x^2)$ and the above polynomial factorization.

(iii) If $2-a_2$ is a square, then let $a_3 \in \fq$ satisfy $a_3^2 = 2-a_2$. Because $a_2^2 = -5$ and $q>3$, we have $a_3^{-2} = \frac{1}{9} (2+a_2)$. Hence $2-9a_3^{-2} = -a_2$. Therefore we have
\begin{eqnarray*}
&&\left( x^4 + a_3 x^3 +  x^2 + 3a_3^{-1}  x +1 \right)
\left( x^4 - a_3 x^3 + x^2 -3 a_3^{-1}  x +1 \right) \\
&=& x^8 +(2-a_3^2) x^6 + (3-6) x^4 + (2-9a_3^{-2})x^2 + 1 \\
&=& x^8 + a_2 x^6 + 3\rho_1 x^4 + a_2 \rho_1 x^2 + 1  
\end{eqnarray*}
Similarly, if $2-a_2$ is a non-square, then $2+a_2$ is also non-square and thus $-(2+a_2)$ is a square. Let $a_3^2 = -(2+a_2)$. Then $a_3^{-2} = -\frac{1}{9} (2-a_2)$ as $a_2^2 = -5$ and $q >3$. Hence $2+9a_3^{-2} = a_2$. 
 Therefore we have
\begin{eqnarray*}
&&\left( x^4 + a_3 x^3 -  x^2 + 3a_3^{-1}  x +1 \right)
\left( x^4 - a_3 x^3 - x^2 -3 a_3^{-1}  x +1 \right) \\
&=& x^8 +(-2-a_3^2) x^6 + (3-6) x^4 + (-2-9a_3^{-2})x^2 + 1 \\
&=& x^8 + a_2 x^6 + 3\rho_1 x^4 + a_2 \rho_1 x^2 + 1  
\end{eqnarray*}

Because $Q_{40}(x)$ is factorized into monic irreducible quartic polynomials under our assumption, we are done.

(iv)
Consider the irreducible factorization of the following format
\begin{eqnarray*}
&& x^8 + a_3 x^6 + b_3 x^4 + c_3 x^2 + 1  \\
&=&\left( x^4 + d_3 x^3 + d_2 x^2 + d_1  x +d_0 \right)
\left(   x^4 + e_3 x^3 + e_2 x^2 + e_1  x +e_0 \right) 
\end{eqnarray*}

Because $d_0$ and $e_0$ are of form $\beta^{1+q+q^2+q^3}$ for some primitive $80$-th root of unity $\beta$ under our assumption, $d_0 = e_0 = \pm 1$. Also one can easily show that $e_3 = -d_3$ and $e_1 = -d_1$ as coefficients of $x^7$ and $x$ vanish in the product. This forces that $d_2 = e_2$. This means that the factorization of $x^8 + a_3 x^6 + b_3 x^4 + c_3 x^2 + 1$ (let $c_3 = 3 a_3^{-1}$) can only be one of the following two ways (either $e_0=d_0 =1$ or $e_0=d_0 = -1$). 
\begin{eqnarray*}
&& x^8 + a_3 x^6 + b_3 x^4 + c_3 x^2 + 1  \\
&=&\left( x^4 + a_4 x^3 + b_4 x^2 + c_4  x \pm 1 \right)
\left(   x^4 -a_4 x^3 + b_4 x^2 -c_4  x  \pm 1 \right) 
\end{eqnarray*}

Comparing coefficients of $x^6, x^4, x^2$ on both sides we have 
\begin{center}
$
\begin{array}{lll}
2b_4 -a_4^2 &=& a_3\\
2 +b_4^2 - 2a_4 c_4  &= & b_3\\
2b_4 -c_4^2 &=& c_3
\end{array}
$
or 
$
\begin{array}{lll}
2b_4 -a_4^2 &=& a_3\\
-2 +b_4^2  - 2a_4c_4  &= & b_3\\
-2b_4 -c_4^2 &=& c_3
\end{array}
$
\end{center}

First we consider the case that $2-a_2$ is a square. In this case, $b_3 =1$.   We now show that if $2-a_2$ is a square then $-2$ is also a square. Indeed, let $a_3^2 = 2-a_2$. Then $a_3^{-2} = \frac{1}{9} (2+a_2)$ and thus $(a_3 -3a_3^{-1})^2 = a_3^2 -6 + 9a_3^{-2} = -2$. Hence $-2$ is a square. Note that the Legendre symbol $\left( -2 \atop p \right) = 1$ iff $p \equiv 1, 3 \pmod{8}$. Let $q=20k+3$. So $k$ is even in this case. Let $\alpha^2 = -2$.

Let us first consider 
 \[
\begin{array}{lll}
2b_4 -a_4^2 &=& a_3\\
-2 +b_4^2 - 2a_4 c_4  &= & 1\\
- 2b_4 -c_4^2 &=& c_3
\end{array}
\]
This case corresponds to $d_0 = e_0 = -1$.  Since $a_4^2 = 2b_4 - a_3$ and $c_4^2 = -2b_4 -c_3$, we obtain
$a_4^2 c_4^2 = (2b_4 -a_3)(-2b_4 -c_3)$. Using $4a_4^2c_4^2 = (b_4^2-3)^2$, we can obtain a quartic equation on $b_4$ only. That is,
\begin{equation}
\label{eqn1}
 b_4^4 + 10 b_4^2 - 8(a_3 - c_3) b_4 + 9-4a_3c_3 =0 
\end{equation}

Note $a_3 c_3 = 3$. 
Since $(a_3 -c_3)^2 = a_3^2 - 2a_3 c_3 + c_3^2 = 2-a_2 - 6 + \frac{1}{9} a_3^{-2} = 2-a_2 -6 + 2 + a_2 = -2$, Equation~(\ref{eqn1}) reduces $b_4^2 + 10b_4^2 -8\alpha b_4 -3 =0$. Moreover,
$b_4^4 + 10b_4^2 - 8 \alpha b_4 -3 = (b_4^2 - 2\alpha b_4 -1)(b_4^2  + 2\alpha b_4 +3) =0$.  Since $-1$ is a non-square in $\fq$,  we must have
$b_4^2 + 2 \alpha b_4 + 3 =0$ and $b_4$ must be one of $\alpha \pm a_2$. 

Secondly we  consider 
 \[
\begin{array}{lll}
2b_4 -a_4^2 &=& a_3\\
2 +b_4^2 - 2a_4 c_4  &= & 1\\
 2b_4 -c_4^2 &=& c_3
\end{array}
\]

Since $a_4^2 = 2b_4 - a_3$ and $c_4^2 = 2b_4 -c_3$, we obtain
$a_4^2 c_4^2 = (2b_4 -a_3)(2b_4 -c_3)$. Using $4a_4^2c_4^2 = (b_4^2 +1)^2$, we can obtain a quartic equation on $b_4$ only. That is,
\begin{equation}
\label{eqn2}
 b_4^4 -14 b_4^2 + 8(a_3 +c_3) b_4 + 1 -4a_3c_3 =0 
\end{equation}

It is easy to check (for example, MAPLE) that Equation~(\ref{eqn2}) has no roots in $\fq$. Since $Q_{80}(x)$ must factor into a product of quartic irreducible polynomials and we know one of above two cases holds,  this shows that we must be able to  find $a_4, b_4, c_4 \in \fq$ such that $a_4^2 = 2b_4 -a_3$, $c_4 = (b_4^2 -3)(2a_4)^{-1}$, $c_4^2 = -2b_4-c_3$,  where $b_4$ must be one of $\alpha \pm a_2 \in \fq$.  

Similarly, if $2-a_2$ is non-square, then we can show that $2$ is a square element. Let $\beta^2 = 2$. In this case, $(a_3 -c_3)^2 = -10$.  If the factorization holds for $d_0 = e_0 = -1$, we have the corresponding equation
\begin{equation}
\label{eqn3}
 b_4^4 + 14 b_4^2 - 8(a_3 - c_3) b_4 + 1-4a_3c_3 =0.
\end{equation}
Otherwise, $(a_3+c_3)^2 = 2$ and we have the equation
\begin{equation}
\label{eqn4}
 b_4^4 -10 b_4^2 + 8(a_3 +c_3) b_4 + 9 -4a_3c_3 =0.
\end{equation}

In fact, it is easy to check that Equation~(\ref{eqn3}) has no solution in $\fq$, but Equation~(\ref{eqn4}) simplifies to $b_4^4 - 10 b_4^2 + 8\beta b_4-3 = (b_4^2 + 2 \beta b_4 - 3)(b_4 - \beta +1)(b_4 - (\beta+1)) =0$. Because $5$ is  a non-square element in $\fq$,  it follows that $b_4^2 + 2 \beta b_4 - 3$ is irreducible and thus   $b_4$ must be one of $\beta \pm 1$ where $\beta^2 = 2$. Again, we must be able to  find $a_4, b_4, c_4 \in \fq$ such that $a_4^2 = 2b_4 -a_3$, $c_4 = (b_4^2 +3)(2a_4)^{-1}$, and $c_4^2 = 2b_4 -c_3$  satisfying Equation~(\ref{eqn4}).

(v) Each quartic factor in $Q_{80}(x)$ has degree $4$ and order $80$, by Theorem~\ref{general} and $Q_{2^n 5}(x) = Q_{80}(x^{2^{n-4}})$, hence the proof is completed.
\end{proof}

We remark that from the above proof,  if $q =20k+3$ where $k>0$ then $k$ is even when $2-a_2$ is a square and $k$ is odd when $-2-a_2$ is a square. 

\begin{corol} \label{irred-mod-3}
Let $q=20k+3$ where $k>0$ and $a_{2}^2 = -5$.

(i)  If $k$ is even, then the polynomial 
\[
x^{2^{n-2}} + a_4 x^{3\cdot 2^{n-4}} + b_4 x^{2^{n-3}} + c_4 x^{2^{n-4}} - 1
\]
is irreducible over $\fq$ for any $n \geq 4$, where $a_4, b_4, c_4$ satisfy (\ref{cond4a}) for any $a_3^2 = 2-a_2$ and $c_3 = 3a_3^{-1}$.

(i)  If $k$ is odd, then the polynomial 
\[
x^{2^{n-2}} + a_4 x^{3\cdot 2^{n-4}} + b_4 x^{2^{n-3}} + (b_4^2+3)(2a_4)^{-1}  x^{2^{n-4}} + 1
\]
is irreducible over $\fq$ for any $n \geq 4$, where $a_4, b_4, c_4$ satisfy (\ref{cond4b}) for any $a_3^2 = -2-a_2$ and $c_3 = 3a_3^{-1}$.  
\end{corol}

\begin{theorem}\label{factor-3}
Let $q=3$. We have the following

(i) If $0 \leq n \leq 1$, then  
\[
Q_{2^n 5} (x) = \prod_{\rho_n \in \Omega(2^n)}  \left( x^4 + \rho_n x^3 + \rho_n^2 x^2 + \rho_n^3 x + \rho_n^4 \right).
\]

(ii) 
\[
Q_{2^2 5} (x) = \prod_{a_2^2 =1}  \left( x^4 + a_2 x^3 -a_2 x +1 \right).
\]

(iii) 
\[
Q_{2^3 5} (x) = \prod_{a_3^2 =1}  \left( x^4 + a_3 x^3 + x^2 +1 \right) \left( x^4 + x^2 + a_3 x +1 \right).
\]

(iv) 
\[
Q_{2^4 5} (x) = \prod_{a_4^2 =1}  \left( x^4 + a_4 x^3 + 2  \right) \left( x^4 + a_4 x +2 \right)  \left( x^4 + a_4 x^3  + x^2 -a_4 x +2 \right) \left( x^4 +a_4 x^3 - x^2 - a_4 x +2 \right).
\]

(v) If $n > 4$, then 
\begin{eqnarray*}
Q_{2^n 5} (x) = & \prod_{a_4^2 =1}  \left( x^{2^{n-2}} + a_4 x^{3\cdot 2^{n-4}} + 2  \right)  \left( x^{2^{n-2}} + a_4 x^{3\cdot 2^{n-4}}  + x^{2^{n-3}} -a_4 x^{2^{n-4}} +2 \right) \\
& \left( x^{2^{n-2}} + a_4 x^{2^{n-4}} +2 \right) \left( x^{2^{n-2}} + a_4 x^{3\cdot 2^{n-4}}  - x^{2^{n-3}} -a_4 x^{2^{n-4}} +2 \right). \\
\end{eqnarray*}
\end{theorem}

\begin{proof}
It is easy to check  directly by computer.
\end{proof}
\begin{corol} \label{irred-3}
For any $n >4$ and $a = \pm 1$,  we have the following families of irreducible polynomials  of degree $2^{n-2}$ over $\F_3$. 

(i) $x^{2^{n-2}} + a x^{3\cdot 2^{n-4}} + 2$;

(ii) $x^{2^{n-2}} + a x^{2^{n-4}} +2$;

(iii) $x^{2^{n-2}} + a x^{3\cdot 2^{n-4}}  + x^{2^{n-3}} -a x^{2^{n-4}} +2$;

(iv) $x^{2^{n-2}} + a x^{3\cdot 2^{n-4}}  - x^{2^{n-3}} -a x^{2^{n-4}} +2 $. 
\end{corol}

\subsection{Case $q \equiv 7 \pmod{20}$}
\label{new3}

Let $q = 20k+7$. If $k$ is odd, then $L_1=1$, $L_2 = L_1+2 = 3$, and $L_4=L_2 + 1 = 4$. The factorization of $Q_{2^n5}(x)$ over $\fq$ behaves the same as that over $\F_{20k+3}$  where $k>0$. If $k$ is  even, then $L_1=1$, $L_2 = L_1+ 3 = 4$, and $L_4=L_2 + 1 = 5$. The factorization of $Q_{2^n5}(x)$ over $\fq$ behaves almost the same as that over $\F_{20t+3}$  where $t>0$ except that the only difference happens when $k$ is even and $n=5$ because $2$ is a square in this case. Again, $Q_{2^5 5}(x) = Q_{2^4 5}(x^2)$,  we need to factor
\[
\displaystyle{ \prod_{a_{2}^2 = -5} \prod_{a_3^2 = - 2  - a_{2} \atop{c_3 = 3a_3^{-1}}}  
\prod_{ a_4, b_4, c_4} }
\left( x^8 + a_4 x^6 + b_4 x^4 + c_4  x^2 + 1 \right)
\]
over $\fq$ where $q=20k+7$ and $k$ is even, where $a_4, b_4, c_4$ satisfy (\ref{cond4b}). Hence the factorization of $Q_{2^5 5}(x)$ is either
\begin{equation}\label{eqn5}
\displaystyle{ \prod_{a_{2}^2 = -5} \prod_{a_3^2 = - 2  - a_{2} \atop{c_3 = 3a_3^{-1}}}  
\prod_{a_4, b_4, c_4}  \prod_{a_5, b_5, c_5}}
\left( x^4 + a_5 x^3 + b_5 x^2 +  c_5 x^2 + 1 \right),
\end{equation}
where $a_5, b_5, c_5$ satisfies $a_5^2 = 2b_5-a_4$ is a square, $c_5 = 2b_5-c_4$ is a square, $c_5 = (b_5^2 +2-b_4)(2a_5)^{-1}$, and
\[
b_5^4 + (-12-2b_4) b_5^2 + 8(a_4+c_4) b_5 + (2 - b_4)^2 -4a_4c_4 =0, 
\]
or
\begin{equation}\label{eqn6}
\displaystyle{ \prod_{a_{2}^2 = -5} \prod_{a_3^2 = - 2  - a_{2} \atop{c_3 = 3a_3^{-1}}}  
\prod_{a_4, b_4, c_4}  \prod_{a_5, b_5, c_5} }
\left( x^4 + a_5 x^3 + b_5 x^2 + c_5 x^2 - 1 \right),
\end{equation}
where $a_5, b_5, c_5$ satisfies  $a_5^2 = 2b_5-a_4$ is a square,  $c_5 = -2b_5-c_4$ is a square, $c_5 = (b_5^2 -2-b_4)(2a_5)^{-1}$, and
\[
b_5^4 + (12-2b_4) b_5^2 - 8(a_4-c_4) b_5 + (2+b_4)^2 -4a_4c_4 =0. 
\]

In these cases, the expressions for $a_5$, $b_5$ are more complicated than those of $a_4$ and $b_4$ and we are not trying to write them explicitly. Similarly, when $q = 20k+7$, we have that  $k$ is odd if $2-a_2$ is a square and $k$ is even if $-2-a_2$ is a square. 

We summarize the result as follows:

\begin{theorem}\label{factor-mod-7}
Let $q\equiv 7 \pmod{20}$.  Then we have the following factorization of cyclotomic polynomials $Q_{2^n 5} (x)$  over $\fq$.

(i) If $0 \leq n \leq 1$, then  
\[
Q_{2^n 5} (x) = \prod_{\rho_n \in \Omega(2^n)}  \left( x^4 + \rho_n x^3 + \rho_n^2 x^2 + \rho_n^3 x + \rho_n^4\right).
\]

(ii) If $n = 2$, then 
\[
Q_{2^n 5} (x) =  \prod_{a_{n}^2 = -5 }  \left( x^4 + a_n x^3 + 3\rho_{1} x^2 + a_n \rho_1 x + 1\right).
\]

(ii) If $n = 3$, then 
\[
Q_{2^n 5} (x) =
 \left\{ 
\begin{array}{ll} 
 \displaystyle{ \prod_{a_{2}^2 = -5} \prod_{ a_3^2 = 2 b_3 - a_{2} \atop{b_3= 1} } } \left( x^4 + a_3 x^3 + b_3 x^2 + 3a_3^{-1}  x + 1 \right), &
if~ 2  - a_{2} ~is~ a ~square;  \\
& \\
\displaystyle{ \prod_{a_{2}^2 = -5} \prod_{a_3^2 = 2 b_3 - a_{2} \atop{b_3= - 1} } } \left( x^4 + a_3 x^3 + b_3 x^2 + 3a_3^{-1}  x + 1 \right), &
 if~ -2 - a_{2} ~is~a~ square. 
\end{array}
\right.
\]

(iii) If $ n = 4$, then 

\[
Q_{2^n 5} (x) = \left\{ 
\begin{array}{ll} 
 \displaystyle{ \prod_{a_{2}^2 = -5} \prod_{ a_3^2 = 2  - a_{2} \atop{c_3 = 3a_3^{-1}} } \prod_{a_4, b_4, c_4 } } \left( x^4 + a_4 x^3 + b_4 x^2 + c_4  x - 1 \right), & 
if~ 2 - a_{2} ~is~ a~square; \\
& \\
\displaystyle{ \prod_{a_{2}^2 = -5} \prod_{a_3^2 = - 2  - a_{2} \atop{c_3 = 3a_3^{-1}}}  
\prod_{a_4, b_4, c_4} }
\left( x^4 + a_4 x^3 + b_4 x^2 + c_4 x + 1 \right), &
 if~ - 2 - a_{2} ~is~a~ square. 
\end{array}
\right.
\]
where $a_4, b_4, c_4$ satisfy either (\ref{cond4a}) if $2-a_2$ is a square,  or (\ref{cond4b}) if $-2-a_2$ is a square. 

(iv) If $n = 5$, then 

\[
Q_{2^n 5} (x) = \left\{ 
\begin{array}{ll} 
 \displaystyle{ \prod_{a_{2}^2 = -5} \prod_{ a_3^2 = 2  - a_{2} \atop{c_3 = 3a_3^{-1}} } \prod_{a_4, b_4, c_4} } \left( x^8 + a_4 x^6 + b_4 x^4 + c_4 x^2 - 1 \right), & 
if~ 2 - a_{2} ~is~ a~ square; \\
& \\
\displaystyle{ \prod_{a_{2}^2 = -5} \prod_{a_3^2 = - 2  - a_{2} \atop{c_3 = 3a_3^{-1}}}  
\prod_{a_4, b_4, c_4} \prod_{a_5, b_5, c_5}}
\left( x^4 + a_5 x^3 + b_5 x^2 + c_5 x \pm 1 \right), &
 if~ - 2 - a_{2} ~is~ a~square,
\end{array}
\right.
\]
where $a_4, b_4, c_4$ satisfy either (\ref{cond4a}) if $2-a_2$ is a square,  or (\ref{cond4b}) if $-2-a_2$ is a square, and  $a_5, b_5, c_5$ are given in (\ref{eqn5}) or (\ref{eqn6}).

(iv) If $n > 5$, then $Q_{2^n 5} (x)$ can be factorized as 

\[
\left\{ 
\begin{array}{ll} 
 \displaystyle{ \prod_{a_{2}^2 = -5} \prod_{ a_3^2 = 2  - a_{2} \atop{c_3 = 3a_3^{-1}} } \prod_{a_4,b_4, c_4} } \left( x^{2^{n-2}} + a_4 x^{3\cdot 2^{n-4}} + b_4 x^{2^{n-3}} + c_4 x^{2^{n-4}} - 1 \right), & 
if~ 2 - a_{2} ~is~ a~ square; \\
& \\
\displaystyle{ \prod_{a_{2}^2 = -5} \prod_{a_3^2 = - 2  - a_{2} \atop{c_3 = 3a_3^{-1}}}  
\prod_{a_4, b_4, c_4} \prod_{a_5, b_5, c_5}}
\left( x^{2^{n-3}} + a_5 x^{3\cdot 2^{n-5}} + b_5 x^{2^{n-4}} + c_5  x^{2^{n-5}} \pm 1 \right), &
 if~ - 2 - a_{2} ~is~a ~ square,
\end{array}
\right.
\]
where $a_4, b_4, c_4$ satisfy either (\ref{cond4a}) if $2-a_2$ is a square,  or (\ref{cond4b}) if $-2-a_2$ is a square, and $a_5, b_5, c_5$ are given in (\ref{eqn5}) or (\ref{eqn6}).
\end{theorem}

\section{Conclusion}
In this paper, we obtain the explicit factorization of cyclotomic polynomials of $Q_{2^n r}(x)$ over finite fields 
where $r=5$ and construct several classes of irreducible polynomials of degree $2^{n-2}$ with fewer than $5$ terms. 
Our approach is recursive,  i.e.,  we derive the factorization of $Q_{2^k r}(x)$ from  the factorization 
of $Q_{2^{k-1} r}(x^2)$. We show that we can do it with at most $L_{\phi(r)} = v_2(q^{\phi(r)} -1)$ iterations. 
A key component of our approach for $r=5$ is to factor certain types of polynomials of degree $8$ into two quartic 
irreducible polynomials.  It would be more desirable to obtain explicit factors of $Q_{2^n r}(x)$  for arbitrary $r$.  One would expect that it  involves 
the factorization of certain types of polynomials of degree $2m$ where $m \mid \phi(r)$ into a product irreducible  polynomials of degree less than or equal to $m$.   Another contribution of this paper is the construction of several classes of  irreducible polynomials over finite fields with at most $5$ nonzero terms. The reciprocals of these irreducible polynomials are also of format 
$x^{2^{n-2}} + g(x)$ such that the degree of $g$ is at most $4$, which could have potential applications as mentioned 
in \cite{GaoHowellPanario}. Finally we note that  one can also construct more classes of irreducible polynomials 
for other choices of $r$  as a consequence of Theorem~\ref{general}.

\section*{Acknowledgement}
Part of the work was done while Qiang Wang took his sabbatical leave at the Center for Advanced Study, Tsinghua University. He wants to thank the Center for their warm hospitality.

\end{document}